\documentclass{amsart}
\usepackage{amssymb}

\usepackage{a4}
\begin{document}
\newtheorem{theorem}{Theorem}[section]
\newtheorem{definition}[theorem]{Definition}
\newtheorem{lemma}[theorem]{Lemma}
\newtheorem{example}[theorem]{Example}
\newtheorem{remark}[theorem]{Remark}
\makeatletter
 \renewcommand{\theequation}{%
 \thesection.\alph{equation}}
 \@addtoreset{equation}{section}
 \makeatother
\title[Affine curvature homogeneous 3-dimensional Lorentz Manifolds]
{Affine curvature homogeneous 3-dimensional Lorentz Manifolds}
\author{P. Gilkey and S. Nik\v cevi\'c}
\begin{address}{PG: Mathematics Department, University of Oregon,
Eugene Or 97403 USA.\newline Email: {\it gilkey@darkwing.uoregon.edu}}
\end{address}
\begin{address}{SN: Mathematical Institute, SANU,
Knez Mihailova 35, p.p. 367,
11001 Belgrade,
Serbia and Montenegro.
\newline Email: {\it stanan@mi.sanu.ac.yu}}\end{address}
\begin{abstract} We study a family of $3$-dimensional Lorentz manifolds. Some members of the family
are $0$-curvature homogeneous, $1$-affine curvature homogeneous, but not $1$-curvature homogeneous. Some are
$1$-curvature homogeneous but not $2$-curvature homogeneous. All are $0$-modeled on indecomposible
local symmetric spaces. Some of the members of the family are geodesically complete, others are not. All have
vanishing scalar invariants.\end{abstract}
\keywords{{$0$-curvature homogeneous, $1$-curvature homogeneous, $1$-affine curvature homogeneous, $3$-dimensional
Lorentz manifold, VSI manifold}, geodesic completeness.
\newline 2000 {\it Mathematics Subject Classification.} 53B20}
\maketitle
\section{Introduction}\label{sect-1}

\subsection{Affine manifolds} We say that $\mathcal{A}:=(M,\nabla)$ is an {\it affine manifold} if $\nabla$ is a torsion free
connection on the tangent bundle
$TM$ of a smooth $m$-dimensional manifold $M$. Let 
$$\mathcal{R}(X,Y):=\nabla_X\nabla_Y-\nabla_Y\nabla_X-\nabla_{[X,Y]}$$ 
be the associated
curvature operator. We say that
$\mathcal{A}$ is {\it locally affine homogeneous} if given any points $P,Q\in M$, there is a diffeomorphism $\Phi_{P,Q}$ from
a neighborhood of $P$ to a neighborhood of $Q$ with $\Phi_{P,Q}(P)=Q$ so that $\Phi_{P,Q}^*\nabla=\nabla$. We say that
$\mathcal{A}$ is {\it locally  $k$-affine curvature homogeneous} if given any points $P,Q\in M$, there is a linear isomorphism
$\phi_{P,Q}$ from
$T_PM$ to $T_QM$ so that $\phi_{P,Q}^*\nabla^i\mathcal{R}_Q=\nabla^i\mathcal{R}_P$ for $0\le i\le k$. By taking
$\phi_{P,Q}=(\Phi_{P,Q})_*$, it is clear that any locally affine homogeneous manifold is locally $k$-affine curvature homogeneous
for all $k$. What is perhaps somewhat surprising is that given $k$, there exists a $k$-affine curvature homogeneous
manifold $\mathcal{A}_k$ of dimension $a(k)$ which is not locally affine homogeneous, see, for example, the discussion in
\cite{GNpre,KO99,O96,O97}; one has that $a(k)\rightarrow\infty$ as $k\rightarrow\infty$.

\subsection{Pseudo-Riemannian manifolds} There are similar notions in the metric context. Let $\mathcal{M}:=(M,g)$ be a
pseudo-Riemannian manifold of signature
$(p,q)$. We take $\nabla$ to be the Levi-Civita connection and let $R\in\otimes^4T^*M$ be the associated curvature tensor: 
$$R(X,Y,Z,W):=g(\mathcal{R}(X,Y)Z,W)\,.$$
We say that $\mathcal{M}$ is {\it locally homogeneous} if given any points $P,Q\in M$, there is an isometry
$\Phi_{P,Q}$ from a neighborhood of $P$ to a  neighborhood of $Q$ with $\Phi_{P,Q}(P)=Q$. We say
that $\mathcal{M}$ is $k$-curvature homogeneous if given any two points
$P,Q\in M$, there is an isometry $\phi_{P,Q}$ from $T_PM$ to $T_QM$ so that $\phi_{P,Q}^*\nabla^iR_Q=\nabla^iR_P$ for $0\le
i\le k$. As $(g,R)$ determines $\mathcal{R}$, locally  homogeneous (resp. $k$-curvature homogeneous) manifolds
are locally affine homogeneous (resp. $k$-affine curvature homogeneous). We refer to the discussion in \cite{BVK96} for a
review of some of the literature in this subject.

Given $k$, there is a pseudo-Riemannian manifold $\mathcal{M}_k$ of dimension $m(k)$ which is $k$-curvature homogeneous (and hence
$k$-affine curvature homogeneous) but not locally affine homogeneous (and hence not locally homogeneous) \cite{GNpre}; one has that
$m(k)\rightarrow\infty$ as $k\rightarrow\infty$. 

If,  however, $m$ is bounded, one has the following result due
to Singer \cite{S60} in the Riemannian ($p=0$) setting and to Podesta and Spiro
\cite{PS04} in the general setting:
\begin{theorem}\label{thm-1.1} There exists an integer $k_{p,q}$ so that if $\mathcal{M}$ is a geodesically complete simply
connected pseudo-Riemannian manifold of signature $(p,q)$ which is $k_{p,q}$-curvature homogeneous, then $\mathcal{M}$ is
homogeneous.
\end{theorem}

We refer to Opozoda \cite{O97} for a similar result in the affine setting; there is an additional technical hypothesis which must
be imposed.

\subsection{Vanishing scalar invariants} Adopt the Einstein convention and sum over repeated indices. We can construct scalar
invariants by contracting indices. For example, the scalar curvature $\tau$, the norm $|\rho|^2$ of the Ricci tensor, and the norm
$|R|^2$ of the full curvature tensor
 are scalar invariants defined by:
\begin{eqnarray*}
&&\tau:=g^{i_1j_1}g^{i_2j_2}R_{i_1i_2j_2j_1},\\
&&|\rho|^2:=g^{i_1j_1}g^{i_2j_2}g^{i_3j_3}g^{i_4j_4}R_{i_1i_2j_2i_3}R_{j_1i_4j_4j_3},\\
&&|R|^2:=g^{i_1j_1}g^{i_2j_2}g^{i_3j_3}g^{i_4j_4}R_{i_1i_2i_3i_4}R_{j_1j_2j_3j_4}\,.
\end{eqnarray*}
By Weyl's theorem \cite{W46}, all universal polynomial scalar invariants of the covariant derivatives of the curvature tensor
arise in this way; thus such invariants are called {\it Weyl scalar invariants}. We say that a pseudo-Riemannian manifold is {\it
VSI} if all the scalar Weyl invariants vanish. This is not possible for non-flat manifolds in the Riemannian setting but is
possible in the higher signature setting, see, for example, the discussion in
\cite{KM96,PPCM02}.

\subsection{Lorentzian manifolds} In this paper, we shall deal with the $3$-dimensional Lorentzian setting -- i.e. signature
$(1,2)$. We shall be discussing a number of tensors. For the sake of brevity, we shall only give the non-zero components up to the
usual symmetries. 
Let $\{x,y,\tilde
x\}$ be coordinates on
$\mathbb{R}^3$. Let $f=f(y)$ be a smooth function on $\mathbb{R}$ and let $\mathcal{M}_{f}:=(\mathbb{R}^3,g_{f})$ where
$g_{f}$ is the Lorentz metric on
$\mathbb{R}^3$ given by:
$$
g_f(\partial_x,\partial_x)=-2f(y)\quad\text{and}\quad g_f(\partial_x,\partial_{\tilde x})=g_f(\partial_y,\partial_y)=1\,.
$$

Let $\mathcal{S}_\varepsilon$ be defined by $f_\varepsilon(y):=\frac12\varepsilon y^2$ for $\varepsilon=\pm1$. 

\begin{theorem}\ \label{thm-1.2}
\begin{enumerate}
\smallbreak\item All scalar Weyl invariants of $\mathcal{M}_f$ vanish. 
\item 
 $\mathcal{S}_\varepsilon$ is an indecomposable local symmetric space.
\item If $f^{\prime\prime}(y)\ne0$, then $\mathcal{M}_f$ is $0$-curvature modeled on $\mathcal{S}_\varepsilon$ for
$\varepsilon=\operatorname{sign}(f^{\prime\prime})$.
\item Assume that $f^{\prime\prime}(y)\ne0$ and that $f^{\prime\prime\prime}(y)\ne0$ for all
$y\in\mathbb{R}$.\begin{enumerate}
\item $\mathcal{M}_f$ is $1$-affine curvature homogeneous.
\item $\mathcal{M}_f$ is $1$-curvature homogeneous if and only if $f^{\prime\prime}=ae^{by}$.
\item The following assertions are equivalent:
\begin{enumerate}
\item $\mathcal{M}_f$ is locally homogeneous.
\item $\mathcal{M}_f$ is $k$-curvature homogeneous for all $k$.
\item $\mathcal{M}_f$ is $2$-curvature homogeneous.
\item $\mathcal{M}_f$ is $2$-affine curvature homogeneous
\item $f^{\prime}=ae^{by}$.
\end{enumerate}
\end{enumerate}\end{enumerate}
\end{theorem}

\subsection{Completeness} Let $\exp_P:T_PM\rightarrow M$ be the exponential map. We say that an affine manifold $\mathcal{A}$ is
{\it geodesically complete} if all geodesics extend for infinite time. 
 
\begin{theorem}\label{thm-1.3}\ 
\begin{enumerate}
\item The manifolds $\mathcal{S}_\pm$ are geodesically complete. 
\item The map $\exp_P$ for $\mathcal{S}_+$ is not surjective for any point
$P\in\mathbb{R}^3$.
\item The map $\exp_P$ for $\mathcal{S}_-$ is a global diffeomorphism from
$T_P\mathbb{R}^3$ to $\mathbb{R}^3$ $\forall P\in\mathbb{R}^3$.
\end{enumerate}\end{theorem}

If $\nabla$ is a torsion free connection, the {\it Jacobi operator} $J_\nabla$ and {\it Ricci form} $\rho_\nabla$ are:
$$J_\nabla(x):\xi\rightarrow\mathcal{R}(\xi,x)x\quad\text{and}\quad
  \rho_\nabla(x,x):=\operatorname{Tr}(J_\nabla(x))\,.$$
An affine manifold $\mathcal{A}$ is said to {\it Ricci explode} if there exists a geodesic $\gamma$ in
$\mathcal{A}$ which is defined for $t\in[0,T)$ where $T<\infty$ so $\lim_{t\rightarrow
T}|\rho(\dot\gamma,\dot\gamma)(t)|=\infty$. Such a manifold is necessarily geodesically incomplete. Furthermore, such a manifold
can not be embedded as an open subset of a geodesically complete affine manifold.

Assume $f^\prime$ never vanishes; by replacing $y$ by $-y$, we may
assume $f^\prime>0$.  The growth of $f^\prime$ at $-\infty$ is crucial.

\begin{theorem}\label{thm-1.4} Assume that $f^\prime(y)>0$ for all $y\in\mathbb{R}$. 
\begin{enumerate}
\item If $\exists$ $C>0$ so $f^\prime(y)\le C|y|$ for $y\le-1$, $\mathcal{M}_f$ is geodesically complete.
\item If $\exists$ $\epsilon,\delta>0$ so $f^{\prime}(y)\ge\epsilon|y|^{1+\delta}$ for
$y\le-1$, $\mathcal{M}_f$ Ricci explodes.\end{enumerate}\end{theorem}

The remainder of this paper is devoted to the proof of these results. In Section \ref{sect-2}, we determine the curvature of the
manifolds $\mathcal{M}_f$ and establish Theorem \ref{thm-1.2}. In Section \ref{sect-3}, we establish Theorem \ref{thm-1.3} by
solving the geodesic equations on $\mathcal{S}_\pm$ quite explicitly. In Section
\ref{sect-4}, we use results from the theory of ordinary differential equations to establish two slightly more general results
from which Theorem
\ref{thm-1.4} will follow.

Various properties of certain of the manifolds in this family have been studied by many authors
\cite{BD00,BV97,CLPT90,KM96,PPCM02}. For example the existence of $1$-curvature homogeneous $3$-dimensional Lorentzian manifolds
which are not locally homogeneous follows from the discussion in \cite{BD00} and the existence of $3$-dimensional VSI Lorentzian
manifolds is established in \cite{PPCM02}. In this paper, we present a unified treatment of a number of results concerning
this family; we discuss some previously known results but also present some new results in affine geometry and deal with questions
of geodesic completeness. We feel this family provides a rich family of examples. In particular, one has:

\begin{example}\label{exm-1.5}\rm For $1\le i\le 3$, let $\mathcal{N}_{i,\pm}:=\mathcal{M}_{f_{i,\pm}}$ where
$$\begin{array}{llll}
f_{1,-}(y)=-e^{-y},&f_{2,-}(y)=-e^{-y}+y,&f_{3,-}(y)=-e^{-y}-e^{-2y}\\
f_{1,+}(y)=e^y,&f_{2,+}(y)=e^y+y,&f_{3,+}(y)=e^y+e^{2y}\,.\end{array}$$
We have $f_i^\prime(y)>0$, $f_i^{\prime\prime}(y)\ne0$, and $f_i^{\prime\prime\prime}\ne0$ for all $y$. We apply the results of
Theorems
\ref{thm-1.1}, \ref{thm-1.2}, \ref{thm-1.3}, and
\ref{thm-1.4} to see:
\begin{enumerate}
\item $\mathcal{S}_-$ is a geodesically complete indecomposible symmetric space.
\item $\mathcal{N}_{1,-}$ is $0$-curvature modeled on $\mathcal{S}_-$, locally
homogeneous, and Ricci explodes.
\item $\mathcal{N}_{2,-}$ is $0$-curvature modeled on $\mathcal{S}_-$, $1$-curvature modeled on $\mathcal{N}_{1,-}$, not
$2$-curvature homogeneous, and  Ricci explodes.
\item $\mathcal{N}_{3,-}$ is $0$-curvature modeled on $\mathcal{S}_-$, not $1$-curvature homogeneous, $1$-affine curvature modeled
on
$\mathcal{N}_1$,  and Ricci explodes.
\item $\mathcal{S}_+$ is a geodesically complete indecomposible symmetric space.
\item $\mathcal{N}_{1,+}$ is $0$-curvature modeled on $\mathcal{S}_+$, geodesically complete, and  homogeneous.
\item $\mathcal{N}_{2,+}$ is $0$-curvature modeled on $\mathcal{S}_+$, $1$-curvature modeled on $\mathcal{N}_{1,+}$, not
$2$-curvature homogeneous,  and geodesically complete.
\item $\mathcal{N}_{3,+}$ is $0$-curvature modeled on $\mathcal{S}_+$, not $1$-curvature homogeneous,
$1$-affine curvature modeled on $\mathcal{N}_{1,+}$,  and geodesically complete.
\end{enumerate}
\end{example}

\section{Curvature}\label{sect-2}

The following Lemma is immediate from the definition:

\begin{lemma}\label{lem-2.1} One has for the manifold $\mathcal{M}_f$ that:
\begin{enumerate}
\item Christoffel symbols:
\begin{enumerate}\item$\nabla_{\partial_x}\partial_x=f^{\prime}\partial_y$.
\item  $\nabla_{\partial_x}\partial_y=\nabla_{\partial_y}\partial_x=-f^{\prime}\partial_{\tilde x}$.\end{enumerate}
\smallbreak\item Components of $R$ and of $\mathcal{R}$:\begin{enumerate}
\item $\mathcal{R}(\partial_x,\partial_y)\partial_y=f^{\prime\prime}\partial_{\tilde x}$.
\item $\mathcal{R}(\partial_x,\partial_y)\partial_x=-f^{\prime\prime}\partial_y$.
\item $R(\partial_x,\partial_y,\partial_y,\partial_x)=f^{\prime\prime}$.\end{enumerate}
\smallbreak\item Components of the Ricci tensor $\rho$:
\begin{enumerate}\item$\rho(\partial_x,\partial_x)=f^{\prime\prime}$.\end{enumerate}
\smallbreak\item Components of $\nabla
R$ and of $\nabla\mathcal{R}$:\begin{enumerate}
\item$\nabla_{\partial_y}\mathcal{R}(\partial_x,\partial_y)\partial_y=f^{\prime\prime\prime}\partial_{\tilde x}$.
\item $\nabla_{\partial_y}\mathcal{R}(\partial_x,\partial_y)\partial_x=-f^{\prime\prime\prime}\partial_y$.
\item $\nabla R(\partial_x,\partial_y,\partial_y,\partial_x;\partial_y)=f^{\prime\prime\prime}$.\end{enumerate}
\smallbreak\item Components of $\nabla^2R$ and of $\nabla^2\mathcal{R}$:\begin{enumerate}
\item 
$\nabla_{\partial y}\nabla_{\partial y}\mathcal{R}(\partial_x,\partial_y)\partial_y=f^{\prime\prime\prime\prime}\partial_{\tilde
x}$.
\item $\nabla_{\partial y}\nabla_{\partial y}\mathcal{R}(\partial_x,\partial_y)\partial_x=
   -f^{\prime\prime\prime\prime}\partial_y$.
\item $\nabla_{\partial x}\nabla_{\partial
x}\mathcal{R}(\partial_x,\partial_y)\partial_y=
   f^{\prime}f^{\prime\prime\prime}\partial_{\tilde x}$.
\item $\nabla_{\partial x}\nabla_{\partial x}\mathcal{R}(\partial_x,\partial_y)\partial_x=
   -f^{\prime}f^{\prime\prime\prime}\partial_y$.
\item $\nabla^2R(\partial_x,\partial_y,\partial_y,\partial_x;\partial_y,\partial_y)=f^{\prime\prime\prime\prime}$.
\item $\nabla^2R(\partial_x,\partial_y,\partial_y,\partial_x;\partial_x,\partial_x)=f^{\prime}f^{\prime\prime\prime}$.
\end{enumerate}\end{enumerate}\end{lemma}

We shall need a technical lemma related to the structure of $\nabla^kR$ when $f$ is a pure exponential. Let
$\nabla^kR(\vec\xi):=\nabla^kR(\xi_1,\xi_2,\xi_3,\xi_4;\xi_5,...,\xi_{4+k})$ for
$\vec\xi=(\xi_1,...,\xi_{4+k})$. We suppose $\xi_i=\partial_x$ or $\xi_i=\partial_y$; there is no need to take
$\xi_i=\partial_{\tilde x}$ since $\nabla^kR(\vec\xi)$ vanishes if any $\xi_i=\partial_{\tilde x}$. Let $\alpha(\vec\xi)$ denote
the number of times that
$\xi_i=\partial_x$.
\begin{lemma}\label{lem-2.2} If $f=ae^{by}$, then
$\nabla^kR(\vec\xi)=\gamma_{\vec\xi}(a,b)e^{\frac12\alpha(\vec\xi)by}$.
\end{lemma}

\begin{proof} We proceed by induction on $k$. Lemma \ref{lem-2.2} follows from Lemma \ref{lem-2.1} when $k=0,1,2$;
$\gamma_{\vec\xi}$ is zero if
$\alpha(\vec\xi)$ is odd. We have by definition that:
\begin{eqnarray}
\nabla^kR(\vec\xi)&=&\xi_{k+4}\nabla^{k-1}R(\xi_1,...,\xi_{k+3})\label{eqn-2.a}\\
&-&\textstyle\sum_{1\le i\le k+3}
  \nabla^{k-1}R(\xi_1,...,\xi_{i-1},\nabla_{\xi_{k+4}}\xi_i,\xi_{i+1},...,\xi_{k+3})\,.\label{eqn-2.b}
\end{eqnarray}

Suppose $\xi_{k+4}=\partial_y$. Let $\vec\eta:=(\xi_1,...,\xi_{k+3})$; $\alpha(\vec\xi)=\alpha(\vec\eta)$. Since
$\nabla_{\partial_y}\xi_i$ is a multiple of
$\partial_{\tilde x}$, the terms in (\ref{eqn-2.b}) vanish and only the term in (\ref{eqn-2.a}) enters. Thus:
\begin{eqnarray*}
&&\nabla^kR(\vec\xi)=\partial_y\nabla^{k-1}R(\vec\eta)
=\partial_y\gamma_{\vec\eta}(a,b)e^{\frac12\alpha(\vec\eta)by}\\
&&\qquad\phantom{.....}=\textstyle\frac12\alpha(\vec\eta)b
\gamma_{\vec\eta}(a,b)e^{\frac12\alpha(\vec\eta)by}
=\gamma_{\vec\xi}(a,b)e^{\frac12\alpha(\vec\xi)by}\quad\text{for}\\
&&\gamma_{\vec\xi}(a,b):=\textstyle\frac12\alpha(\vec\eta)b\gamma_{\vec\eta}(a,b)\,.
\end{eqnarray*}

Suppose that $\xi_{k+4}=\partial_x$. Since $\xi_{k+4}\nabla^{k-1}R(\xi_1,...,\xi_{k+3})=0$, the term in (\ref{eqn-2.a}) vanishes.
Since $\nabla_{\partial_x}\partial_y$ is a multiple of $\partial_{\tilde x}$, we can ignore terms where $\xi_i=\partial_y$. Set 
$$\vec\eta_i:=(\xi_1,...,\xi_{i-1},\partial_y,\xi_{i+1},...,\xi_{k+3})\,.$$
We have
$\nabla_{\partial_x}\partial_x=f^\prime\partial_y$. If $\xi_i=\partial_x$, then $\alpha(\vec\eta_i)=\alpha(\vec\xi)-2$. We compute:
\begin{eqnarray*}
&&\nabla^kR(\vec\xi)
=-abe^{by}\textstyle\sum_{i:\xi_i=\partial_x}\nabla^{k-1}R(\vec\eta_i)
=-abe^{by}\textstyle\sum_{i:\xi_i=\partial_x}\gamma_{\vec\eta_i}(a,b)
e^{\frac12\alpha(\vec\eta_i)by}\\
&&\qquad\phantom{..........}=\gamma_{\vec\xi}(a,b)e^{\frac12\alpha(\vec\xi)by}\quad\text{for}\\
&&\gamma_{\vec\xi}(a,b):=-ab\textstyle\sum_{i:\xi_i=\partial_x}\gamma_{\vec\eta_i}(a,b)\,.
\end{eqnarray*}

The Lemma now follows from these two special cases.
\end{proof}

\begin{proof}[Proof of Theorem \ref{thm-1.2}] We show that all the scalar Weyl invariants of $\mathcal{M}_f$ vanish as follows.
Consider the orthonormal basis
$$e_1^+:=\partial_x+\textstyle\frac12\partial_{\tilde x},\quad e_2^-:=\partial_x-\frac12\partial_{\tilde x},\quad
  e_3^+:=\partial_y\,.$$
We form scalar Weyl invariants by contracting indices in pairs and then summing over repeated indices.
Since $\nabla^kR(...,e_1^+,...)=\nabla^kR(...,e_2^-,...)=\nabla^kR(...,\partial_x,...)$, since $g^{11}=+1$, and
since $g^{22}=-1$, terms where $e_i=e_1^+$ (i.e. $i=1$) and terms where $e_i=e_2^-$ (i.e. $i=2$) appear with opposite signs in any
Weyl summation and cancel; $\nabla^kR(...,e_3^+,...)=0$. Assertion (1) now follows.

If $f$ is quadratic, then $\nabla R=0$. Thus $\mathcal{S}_+$ and
$\mathcal{S}_-$ are local symmetric spaces. The curvature tensor and metric are indecomposible; they are not irreducible
as $\operatorname{Span}\{\partial_{\tilde x}\}$ is invariant  under the isotropy representation. Assertion (2) follows.

Set $\varepsilon_f:=\operatorname{sign}(f^{\prime\prime})$. We say that a basis $\mathcal{B}=\{X,Y,\tilde X\}$ is {\it normalized}
if we have:
\begin{equation}\label{eqn-2.c}
\begin{array}{l}
\text{1) the non-zero components of $g$ are}\quad g(X,\tilde X)=g(Y,Y)=1,\\
\text{2) the non-zero components of $R$ are}\quad R(X,Y,Y,X)=\varepsilon_f,\vphantom{\vrule height 11pt}\\
\text{3) we have}\quad\nabla R(\xi_1,\xi_2,\xi_3,\xi_4;X)=0\ \forall\xi_1,\xi_2,\xi_3,\xi_4\,.\vphantom{\vrule height 11pt}
\end{array}\end{equation}
We may define a normalized basis, and thereby establish Assertion (3), by setting
\begin{equation}\label{eqn-2.d}
X:=|f^{\prime\prime}|^{-1/2}\{\partial_x+f\partial_{\tilde x}\},\quad
  Y:=\partial_y,\quad\tilde X:=|f^{\prime\prime}|^{1/2}\partial_{\tilde x}\,.
\end{equation}

We say $\mathcal{B}$ is {\it affine normalized} if the non-zero components of $\mathcal{R}$ and $\nabla
\mathcal{R}$ are
$$
\begin{array}{ll}
\mathcal{R}(X,Y)Y=\varepsilon_f\tilde X,&\mathcal{R}(X,Y)X=-\varepsilon_fY,\\
\nabla_Y\mathcal{R}(X,Y)Y=\tilde X,&\nabla_Y\mathcal{R}(X,Y)X=-Y\,.
\end{array}
$$
We construct an affine normalized basis by rescaling the coordinate frame. Let $a_1$, $a_2$, and $a_3$ be
constants to be determined. By Lemma \ref{lem-2.1},
\begin{eqnarray*}
&&\mathcal{R}(a_1\partial_x,a_2\partial_y)a_2\partial_y=a_1a_2^2a_3^{-1}f^{\prime\prime}a_3\partial_{\tilde x},\\
&&\mathcal{R}(a_1\partial_x,a_2\partial_y)a_1\partial_x=-a_1^2f^{\prime\prime}a_2\partial_y,\\
&&\nabla_{a_2\partial_y}\mathcal{R}(a_1\partial_x,a_2\partial_y)a_2\partial_y
     =a_1a_2^3a_3^{-1}f^{\prime\prime\prime}a_3\partial_{\tilde x},\\
&&\nabla_{a_2\partial_y}\mathcal{R}(a_1\partial_x,a_2\partial_y)a_1\partial_x=-a_1^2a_2f^{\prime\prime\prime}a_2\partial_y\,.
\end{eqnarray*}
Assume that $f^{\prime\prime}(y)$ and $f^{\prime\prime\prime}(y)$ never vanish. We define an affine normalized basis and prove
Assertion (4a) by setting
\begin{equation}\label{eqn-2.e}
\begin{array}{llll}
X:=a_1\partial_x,&Y:=a_2\partial_y,&\tilde X:=a_3\partial_{\tilde x}&\text{where}\\
a_1:=\{|f^{\prime\prime}|\}^{-1/2},&
a_2:=|f^{\prime\prime}|\{f^{\prime\prime\prime}\}^{-1},&a_3:=a_1a_2^2|f^{\prime\prime}|\,.\vphantom{\vrule height 11pt}
\end{array}\end{equation}
We note for future reference that
\begin{equation}\label{eqn-2.f}
\begin{array}{l}
\nabla_X\nabla_X\mathcal{R}(X,Y)Y=a_1^3a_2^2\nabla_{\partial_x}\nabla_{\partial_x}\mathcal{R}(\partial_x,\partial_y)\partial_y
=a_1^3a_2^2f^\prime f^{\prime\prime\prime}\partial_{\tilde x}\\
\qquad\qquad\qquad\phantom{......}
=f^\prime f^{\prime\prime\prime}\{f^{\prime\prime}\}^{-2}\tilde X\,.\vphantom{\vrule height 11pt}
\end{array}
\end{equation}

We study the relevant symmetry group to construct additional invariants of the $1$-model. Let $\mathcal{B}=\{X,Y,\tilde X\}$
be the normalized basis defined in Equation (\ref{eqn-2.d}). Suppose that
$\mathcal{B}_1=\{X_1,Y_1,\tilde X_1\}$ is another normalized basis. Expand:
\begin{eqnarray*}
&&X_1=a_{11}X+a_{12}Y+a_{13}\tilde X,\\&&Y_1=a_{21}X+a_{22}Y+a_{23}\tilde X,\\
&&\tilde X_1=a_{31}X+a_{32}Y+a_{33}\tilde X\,.
\end{eqnarray*}
Since $R(\xi_1,\xi_2,\xi_3,\tilde X_1)=0$ for any $\xi_1,\xi_2,\xi_3$, we have $a_{31}=a_{32}=0$. Since $\nabla
R(\xi_1,\xi_2,\xi_3,\xi_4;X_1)=0$ for any $\xi_1,\xi_2,\xi_3,\xi_4$, $a_{12}=0$. Thus
$$
X_1=a_{11}X+a_{13}\tilde X,\quad Y_1=a_{21}X+a_{22}Y+a_{23}\tilde X,\quad\tilde X_1=a_{33}\tilde X\,.
$$
As $g(X_1,\tilde X_1)=1$,
$a_{33}a_{11}=1$. As $g(Y_1,\tilde X_1)=0$, $a_{21}=0$.  As $g(X_1,X_1)=0$, $a_{13}=0$. Consequently,
$$
X_1=a_{11}X,\quad Y_1=a_{22}Y+a_{23}\tilde X,\quad\tilde X_1=a_{11}^{-1}\tilde X\,.
$$
As $g(X_1,Y_1)=0$, $a_{23}=0$.
As $g(Y_1,Y_1)=1$, $a_{22}^2=1$. As $R(X_1,Y_1,Y_1,X_1)=\varepsilon_f$, $a_{11}^2a_{22}^2=1$. Thus
$$X_1=a_{11}X,\quad Y_1=a_{22} Y,\quad
\tilde X_1=a_{11}^{-1}\tilde X\quad\text{where}\quad a_{11}^2=a_{22}^2=1\,.
$$
In particular we may use Equation (\ref{eqn-2.d}) and Lemma \ref{lem-2.1} to see:
$$|\nabla R(X_1,Y_1,Y_1,X_1;Y_1)|=|\nabla R(X,Y,Y,X;Y)|=|f^{\prime\prime\prime}\{f^{\prime\prime}\}^{-1}|$$
is an invariant of the $1$-model. This is constant if and only if $f^{\prime\prime\prime}=cf^{\prime\prime}$, i.e.
$f^{\prime\prime}=ae^{by}$. Assertion (4b) now follows.

We now establish Assertion (4c). The following implications are immediate:
$$
   \text{(4c-i)}\Rightarrow\text{(4c-ii)}\Rightarrow\text{(4c-iii)}\Rightarrow\text{(4c-iv)}\,.
$$
Suppose that
$\mathcal{M}_f$ is $2$-affine curvature homogeneous. Let $\mathcal{B}:=\{X,Y,\tilde X\}$ be the affine normalized basis
defined in Equation (\ref{eqn-2.e}). Suppose that
$\mathcal{B}_1:=\{X_1,Y_1,\tilde X_1\}$ is another affine normalized basis. Let
\begin{eqnarray*}
&&\mathcal{I}_0:=\operatorname{Span}_{\xi_1,\xi_2,\xi_3}\{\mathcal{R}(\xi_1,\xi_2)\xi_3\}=\operatorname{Span}\{Y,\tilde X\},\\
&&\mathcal{K}_0:=\{\eta:\mathcal{R}(\eta,\xi_1)\xi_2=0\text{ for all }\xi_1,\xi_2\}=\operatorname{Span}\{\tilde X\},\\
&&\mathcal{K}_1:=\{\eta:\nabla_\eta \mathcal{R}(\xi_1,\xi_2)\xi_3=0\text{ for all }\xi_1,\xi_2,\xi_3\}=
  \operatorname{Span}\{X,\tilde X\}\,.
\end{eqnarray*}
Since these spaces are invariantly defined, we may expand
$$
  X_1=a_{11}X+a_{13}\tilde X,\quad Y_1=a_{22}Y+a_{23}\tilde X,\quad\tilde X_1=a_{33}\tilde X\,.
$$
We have 
$\mathcal{R}(X_1,Y_1)X_1=a_{11}^2\mathcal{R}(X,Y_1)X$ is a multiple of $Y$. Since the basis is normalized, it is also
a multiple of $Y_1$.
Thus $a_{23}=0$. We now compute
$$\begin{array}{ll}
\mathcal{R}(X_1,Y_1)Y_1=a_{11}a_{22}^2a_{33}^{-1}\varepsilon_f\tilde X_1,&
\mathcal{R}(X_1,Y_1)X_1=-a_{11}^2\varepsilon_fY_1,\\
\nabla_{Y_1}\mathcal{R}(X_1,Y_1)Y_1=a_{11}a_{22}^3a_{33}^{-1}\tilde X_1,&
\nabla_{Y_1}\mathcal{R}(X_1,Y_1)X_1=-a_{11}^2a_{22}Y_1\,.\vphantom{\vrule height 11pt}
\end{array}$$
As the basis is normalized, $a_{11}^2=1$, $a_{22}=1$, $a_{33}=a_{11}$. Thus by
Equation (\ref{eqn-2.f}),
\begin{eqnarray*}
&&\nabla_{X_1}\nabla_{X_1}\mathcal{R}(X_1,Y_1)Y_1=a_{11}^3a_{22}^2\nabla_{X}\nabla_{X}\mathcal{R}(X,Y)Y\\
&=&a_{11}^3a_{22}^2
  f^\prime f^{\prime\prime\prime}\{f^{\prime\prime}\}^{-2}\tilde X
 =a_{11}^3a_{22}^2a_{33}^{-1}f^\prime f^{\prime\prime\prime}\{f^{\prime\prime}\}^{-2}\tilde X_1\\
&=&f^\prime f^{\prime\prime\prime}\{f^{\prime\prime}\}^{-2}\tilde X_1\,.
\end{eqnarray*}
This shows that $f^\prime f^{\prime\prime\prime}\{f^{\prime\prime}\}^{-2}$ 
is an invariant of the affine
2-model. Consequently if $\mathcal{M}_{f}$ is $2$-affine curvature homogeneous, then
$$f^{\prime}f^{\prime\prime\prime}\{f^{\prime\prime}\}^{-2}=c\,.$$
If $c=1$, then
$f^{\prime}=ae^{by}$. If $c\ne1$, then $f^\prime=\alpha(y+\beta)^\gamma$ for some $\{\alpha,\beta,\gamma\}$. This later choice is
ruled out as $f^{\prime\prime\prime}$ and $f^{\prime\prime}$ are assumed to be globally defined and non-zero. Thus we may conclude
$f^\prime=ae^{by}$; this establishes the implication:
$$\text{(4c-iv)}\Rightarrow\text{(4c-v)}\,.$$

Finally, we suppose that
$f^\prime=ae^{by}$; we take $f=\frac abe^{by}$.
Consider the normalized basis $\{X,Y,\tilde X\}$ defined in Equation (\ref{eqn-2.d}):
$$X=|ab|^{-1/2}e^{-by/2}\{\partial_x+\textstyle\frac abe^{by}\partial_{\tilde x}\},\quad
  Y=\partial_y,\quad
  \tilde X=|ab|^{1/2}e^{by/2}\partial_{\tilde x}\,.
$$
Let
$\nabla^kR(\vec\eta)=\nabla^kR(\eta_1,\eta_2,\eta_3,\eta_4;\eta_5,...,\eta_{4+k})$ for $\vec\eta=(\eta_1,...,\eta_{4+k})$
where $\eta_i=X$ or $\eta_i=Y$ for $1\le i\le 4+k$. Let $\vec\xi$ be the corresponding string where $X$ and $Y$ are replaced by
$\partial_x$ and $\partial_y$. Let $\alpha(\vec\eta)=\alpha(\vec\xi)$ be the number of times that $X$ or equivalently
$\partial_x$ appear. We apply Lemma \ref{lem-2.2} to see:
\begin{eqnarray*}
\nabla^kR(\vec\eta)&=&|ab|^{-\alpha(\eta)/2}e^{-\alpha(\vec\eta)by/2}\nabla^kR(\vec\xi)\\
&=&|ab|^{-\alpha(\eta)/2}e^{-\alpha(\vec\eta)by/2}\gamma_{\vec\xi}(a,b)e^{\alpha(\vec\xi)by/2}\\
&=&|ab|^{-\alpha(\eta)/2}\gamma_{\vec\xi}(a,b)\,.
\end{eqnarray*}
This shows that $\mathcal{M}_f$ is $k$-curvature homogeneous for all $k$; a local version of Theorem \ref{thm-1.1} now shows
that $\mathcal{M}_f$ is locally homogeneous as desired. Consequently, (4c-v)$\Rightarrow$(4c-i).
\end{proof}

\section{Complete manifolds}\label{sect-3}

Let $\gamma(t)=(x(t),y(t),\tilde x(t))$ be a path in $\mathcal{M}_f$. The geodesic
equation becomes
$$x^{\prime\prime}(t)=0,\quad y^{\prime\prime}(t)=-f^\prime(y(t))x^\prime(t) x^\prime(t),\quad
\tilde x^{\prime\prime}(t)=2f^\prime(y(t))y^\prime(t)x^\prime(t)\,.$$
The first equation yields $x(t)=x_0+x_1t$. The remaining equations then become
$$y^{\prime\prime}(t)=-x_1^2f^{\prime}(y(t)),\quad \tilde x^{\prime\prime}(t)=2x_1f^{\prime}(y(t))y^\prime(t)\,.$$
The equation for $y$ is the crucial one; once $y$ is determined, one can express
\begin{equation}\label{eqn-3.a}
\tilde x(t)=\tilde x_0+t\tilde
x^\prime_0+2x_1\int_{s=0}^t\int_{u=0}^sf^\prime(y(u))y^\prime(u)duds\,.
\end{equation}

\begin{proof}[Proof of Theorem \ref{thm-1.3}] First set $f_+(y)=\frac12y^2$. We then have to solve
$$y^{\prime\prime}(t)=-x_1^2y(t)\,.$$
We show $\mathcal{S}_{+}$ is geodesically complete by solving this equation:
$$
y(t)=\left\{\begin{array}{lll}
y_0+y_0^\prime t&\text{if}&x_1=0,\\
y_0\cos(x_1t)+\frac1{x_1}y_0^\prime\sin(x_1t)&\text{if}&x_1\ne0\,.\vphantom{\vrule height 11pt}
$$
\end{array}\right.$$

A geodesic with $x(0)=x_0$ and $x(1)=x_0+2\pi$ has the form:
$$\gamma(t)=(x_0+2\pi t,y_0\cos(2\pi t)+\textstyle\frac1{2\pi} y_0^\prime\sin(2\pi t),\tilde x(t))\,.$$
Thus $y(1)=y(0)$ and the exponential map is not surjective. This establishes the Assertions of the Lemma concerning
$\mathcal{S}_+$.

Next, we study $f_-(y)=-\frac12y^2$. We then have to solve
$$y^{\prime\prime}(t)=x_1^2y(t)\,.$$
We show $\mathcal{S}_{-}$ is geodesically complete by solving this equation:
$$
y(t)=\left\{\begin{array}{lll}
y_0+y_0^\prime t&\text{if}&x_1=0,\\
\frac12y_0\{e^{x_1t}+e^{-x_1t}\}+\frac1{2x_1}y_0^\prime\{e^{x_1t}-e^{-x_1t}\}&\text{if}&x_1\ne0\,.
\vphantom{\vrule height 11pt}
$$
\end{array}\right.$$

We take $P:=(x_0,y_0,\tilde x_0)$ as the initial point. Suppose $Q:=(x_1,y_1,\tilde x_1)$ is given. The exponential map is given
by setting $t=1$. Thus $x(t)=x_0+t(x_1-x_0)$. If $x_1-x_0=0$, then set $y(t)=y_0+t(y_1-y_0)$. If $x_1-x_0\ne0$, we determine
$y_0^\prime$ uniquely by solving the equation:
$$y_1={\textstyle\frac12}y_0\{e^{x_1-x_0}+e^{x_0-x_1}\}+{\textstyle\frac1{2(x_1-x_0)}}y_0^\prime\{e^{x_1-x_0}-e^{x_0-x_1}\}\,.$$
Once $x$ and $y$ have been determined, we then use Equation (\ref{eqn-3.a}) to solve for $\tilde x_0^\prime$.
This shows that $\mathcal{S}_-$ is geodesically complete and that the exponential map is a diffeomorphism.
\end{proof}

\section{The proof of Theorem \ref{thm-1.4}}\label{sect-4}

Following the discussion in Section \ref{sect-3}, to construct geodesics in the manifold $\mathcal{M}_f$, we must solve the ODE
$$y^{\prime\prime}=-x_1^2f^{\prime}(y)\,.$$

We shall suppose $x_1\ne0$ and set $h=-x_1^2f^\prime$. We begin with:

\begin{lemma}\label{lem-4.1}
Let $h:\mathbb{R}\rightarrow(-\infty,0)$ be smooth. Let $[0,T)$ be the maximal
domain of the solution $y$ to the ODE $y^{\prime\prime}=h(y)$ where $y(0)=y_0$ and $y^\prime(0)=y^\prime_0$. If $T<\infty$,
$$
  \lim_{t\rightarrow T}y(t)=\lim_{t\rightarrow T}y^\prime(t)=-\infty\quad\text{and}\quad
   \limsup_{y\rightarrow T}\frac{h(y(t))}{y(t)}=\infty\,.
$$
\end{lemma}

\begin{proof} Since $y^{\prime\prime}<0$, $y^\prime$ is monotonically decreasing and $y$ is bounded from above
on $[0,T)$. Suppose first that $y$ is bounded from below on $[0,T)$. This implies that $y^{\prime\prime}$ is bounded and hence
$y^\prime$ is bounded as well on
$[0,T)$. Let 
$$
y_1=\lim\inf_{t\rightarrow T}y(t)\quad\text{and}\quad y_1^\prime=\lim_{t\rightarrow T}y^\prime(t)\,.
$$
The fundamental theorem of ODE's shows there exists $\kappa>0$ so that if 
$$
 |z_1-y_1|<\kappa,\quad|z_1^\prime-y_1^\prime|<\kappa,\quad\text{and}\quad s\in(T-\kappa,T)\,
$$
then there exists a solution $z$ to the equation $z^{\prime\prime}=h(z)$ with initial conditions $z(s)=z_1$
and $z^\prime(s)=z_1^\prime$ which is valid on the interval $[s,s+\kappa)$. We choose 
$$
   s\in(T-\textstyle\frac12\kappa,T)\quad\text{so that}\quad
  |y(s)-y_1|<\kappa\quad\text{and}\quad|y^\prime(s)-y_1^\prime|<\kappa\,.
$$
Let $z^{\prime\prime}=h(z)$ be defined on $[s,s+\kappa)$ with $z(s)=y(s)$ and $z^\prime(s)=y^\prime(s)$. Then $z$ extends $y$
to the region
$[0,T+\frac12\kappa)$ which contradicts the assumption that $[0,T)$ was a maximal domain. 

Thus $y$ is not bounded from below on $[0,T)$ so $\lim_{t\rightarrow T}y^\prime(t)=-\infty$. Consequently, $y$ is
monotonically decreasing for $t$ close to $T$ so $\lim_{t\rightarrow T}y(t)=-\infty$ as well. Suppose
$$\limsup_{t\rightarrow T}\frac{h(y(t))}{y(t)}<\infty$$
i.e. that there exists $C<\infty$ so $|h(y(t))|\le C|y(t)|$ on $[t_0,T)$. We then have
$$\left\{\ln|y(t)|\right\}^{\prime\prime}=\left\{\frac{y^\prime(t)}{y(t)}\right\}^\prime=\frac{y^{\prime\prime}(t)}{y(t)}-
\left\{\frac{y^\prime(t)}{y(t)}\right\}^2
=\frac{h(y(t))}{y(t)}-\left\{\frac{y^\prime(t)}{y(t)}\right\}^2\le C\,.$$
This implies $\ln|y(t)|$ is bounded from above and hence $|y(t)|$ is bounded from above on $[t_0,T)$ which is false. This
contradiction shows $\limsup_{t\rightarrow T}\frac{h(y(t))}{y(t)}=\infty$.
\end{proof}

\begin{proof}[Proof of Theorem \ref{thm-1.4} (1)] We suppose that $f^\prime>0$ and that $f^\prime(y)\le C|y|$ for $y\le -1$.
We set $h=-x_1^2f^\prime$. Choose a maximal domain $[0,T)$ for the solution to the ODE $y^{\prime\prime}=h(y)$ with initial
condition $y(0)=y_0$ and
$y^{\prime}(0)=y^\prime_0$. If $T<\infty$, then 
$$\limsup_{t\rightarrow T}\frac{h(y(t))}{y(t)}=\infty$$ 
which is false. Thus $T=\infty$
and $\mathcal{M}_f$ is geodesically complete.\end{proof}

Before proving Theorem \ref{thm-1.4} (2), we must establish:

\begin{lemma}\label{lem-4.2}\ 
\begin{enumerate}
\item Let $\alpha>0$. Let $\{t_n\}_{n\ge1}$ be a sequence of real numbers with $t_1=1$ and with $t_{n+1}-t_{n}\ge n^\alpha$ for
$n\ge1$. Then
$t_n\ge \frac{n^{1+\alpha}}{(1+\alpha)2^{1+\alpha}}$.
\item Let $\epsilon>0$ and $\delta>0$. Suppose that $h(y)<-\epsilon|y|^{1+\delta}$ for $y\le-1$. Let $[0,T)$ be the
maximal domain of definition for the solution $y$ to the ODE $y^{\prime\prime}=h(y)$ with $y(0)=-1$ and
$y^\prime(0)=-1$. Then $T<\infty$ and $\lim_{t\rightarrow T}y(t)=-\infty$.
\end{enumerate}\end{lemma}

\begin{proof} We prove Assertion (1) by induction on $n$; it holds trivially for $n=1$. We take $n\ge2$ and 
use the comparison test to compute:
\begin{eqnarray*}
t_n&>&t_n-t_1=\sum_{k=2}^n\bigg\{t_k-t_{k-1}\bigg\}\ge\int_1^n(x-1)^\alpha dx\\
&=&\frac{(n-1)^{1+\alpha}}{1+\alpha}=
\frac{n^{1+\alpha}}{(   1+\alpha)(1+\frac1{n-1})^{1+\alpha}}\ge\frac{n^{1+\alpha}}{(1+\alpha)2^{1+\alpha}}\,.
\end{eqnarray*}

To prove Assertion (2), we suppose first $T=\infty$ and argue for a contradiction. Choose $\tau$ so that
$$\tau\varepsilon\ge2^{1+\delta/2}(1+\delta/2)\quad\text{and}\quad\tau\ge1\,.$$
With our initial conditions, $y^{\prime\prime}<0$ so $y^\prime$ is monotonically decreasing and $y^\prime\le -1$. This
implies
$y$ decreases monotonically. Let $\Delta_n=\tau\cdot n^{-1-\delta/2}$. Let $s_1=0$ and let $s_{n+1}=s_n+\Delta_n$ for $n\ge2$. As
$\delta>0$,
$$S:=\lim_{n\rightarrow\infty}s_n=\sum_{n=1}^\infty\tau n^{-1-\delta/2}<\infty\,.$$
We wish to show inductively that\begin{enumerate}
\item $y^\prime(s_n)\le-n^{1+\delta/2}$.
\item $y(s_n)\le -n$.
\item $y^\prime(s_{n+1})-y^\prime(s_n)\le -2^{1+\delta/2}(1+{\textstyle\frac12}\delta)n^{\delta/2}$.
\end{enumerate}
The first two statements hold $n=1$ by the choice of our initial conditions. Since $y$ and $y^\prime$ decrease monontonically, we
may estimate
\begin{eqnarray*}
&&y^{\prime\prime}(s)\le-\epsilon|y(s)|^{1+\delta}\le-\epsilon|y(s_n)|^{1+\delta}\le-\epsilon n^{1+\delta}\text{ for
}s\in[s_n,s_{n+1}],\\ &&y^\prime(s_{n+1})-y^\prime(s_n)\le-\Delta_n\epsilon n^{1+\delta}=-\tau
n^{-1-\delta/2}\epsilon n^{1+\delta}\le-2^{1+\delta/2}(1+\delta/2)n^{\delta/2}\,.
\end{eqnarray*}
Thus statements $(1)_n$ and $(2)_n$ imply assertion $(3)_n$.

Statements $(3)_k$ for $1\le k\le n$ together with Assertion (1) imply Statement $(1)_{n+1}$. Finally, 
we use Statement $(1)_n$ together with statement $(2)_n$ to establish Statement $(2)_{n+1}$ by computing:
\begin{eqnarray*}
&&y^\prime(s)\le y^\prime(s_n)\le-n^{1+\delta/2}\quad\text{for}\quad s\ge s_n,\\
&&y(s_{n+1})\le y(s_n)+\Delta_ny^\prime(s_n)\le -n-\tau n^{-1-\delta/2}n^{1+\delta/2}\le-n-1\,.
\end{eqnarray*}
This establishes the truth of all the 3 statements. Thus, $\lim_{s\rightarrow
S}y(s)=-\infty$. This contradicts the assumption that $T=\infty$.

This shows that $y$ must be defined on a maximal domain $[0,T)$ for $T<\infty$; the fact that $\lim_{t\rightarrow T}y(t)=-\infty$
now follows from Lemma \ref{lem-4.1}.
\end{proof}

\begin{proof}[Proof of Theorem \ref{thm-1.4} (2)] Suppose $f^\prime(y)>0$ for all $y$ and that $f^\prime(y)\ge\varepsilon
|y|^{1+\delta}$ for $y\le-1$. Choose a geodesic with
$x(0)=0$, $x^\prime(0)=1$, $y(0)=-1$, and $y^\prime(0)=-1$. We then have the differential equation
$$y^{\prime\prime}=-f^\prime(y)\,.$$
Thus by Lemma \ref{lem-4.2} for some finite time $T$, we have $\lim_{t\rightarrow T}y(t)=-\infty$. Thus $\mathcal{M}_f$ is
geodesically incomplete. We have $\rho(\dot\gamma,\dot\gamma)=f^{\prime\prime}(y(t))$. 

If $|f^{\prime\prime}(y)|\le K$ on $(-\infty,0]$, then $f^\prime(y)\le K|y|+f^\prime(y(0))$ on $(-\infty,0]$ which is false. Thus
$|f^{\prime\prime}(y)|$ is not bounded on $(-\infty,0]$. Since $y(t)\rightarrow-\infty$ as $t\rightarrow T$, 
$f^{\prime\prime}(y(t))$ is not bounded on $[0,T)$. This shows, as desired, that $\mathcal{M}_f$ Ricci explodes. 
\end{proof}

\section*{Acknowledgments} Research of P. Gilkey partially supported by the
Max Planck Institute in the Mathematical Sciences (Leipzig). Research of S. Nik\v cevi\'c partially supported by MM 1646
(Srbija).


\begin{thebibliography}{AAA}
\bibitem{BVK96} E. Boeckx, L. Vanhecke, and O. Kowalski, {\bf Riemannian manifolds of conullity two}, World Scientific (1996).


\bibitem{BD00} P. Bueken and M. Djori\'c, Three-dimensional Lorentz metrics and curvature
   homogeneity of order one, {\it Ann. Global Anal. Geom.} {\bf 18} (2000), 85--103.

\bibitem{BV97} P. Bueken and L. Vanhecke,
Examples of curvature homogeneous Lorentz metrics, {\it Classical Quantum Gravity} {\bf 14} (1997), L93--L96.

\bibitem{CLPT90} M. Cahen, J. Leroy, M. Parker, F. Tricerri, and L. Vanhecke,
Lorentz manifolds modeled on a Lorentz symmetric space, {\it J. Geom. Phys.} {\bf 7} (1990), 571-581.

\bibitem{GNpre} P. Gilkey and S. Nik\v cev\'ic, Complete $k$-curvature homogeneous pseudo-Riemannian manifolds, {to appear \it
Annals Global Analysis and Geometry}, math.DG/0405024.

\bibitem{KO99} O. Kowalski, B. Opozda, and Z. Vl\'{a}\v{s}ek,
Curvature homogeneity of affine connections on two-dimensional manifolds, {\it Colloq. Math.} {\bf 81}
(1999), 123--139.

\bibitem{KM96} A. Koutras and C. McIntosh, A metric with no symmetries or invariants, {\it Classical Quantum Gravity}
{\bf 13} (1996), L47-L49.

\bibitem{O96} B. Opozda, On curvature homogeneous and locally
homogeneous affine connections,
   {\it Proc. Amer. Math. Soc.} {\bf 124} (1996), 1889--1893.

\bibitem{O97} B. Opozda, Affine versions of Singer's theorem on locally
homogeneous spaces, {\it Ann. Global Anal. Geom.} {\bf 15} (1997), 187--199.

\bibitem{PS04} F. Podesta and A. Spiro, Introduzione ai Gruppi di Trasformazioni, 
{\it Volume of the Preprint Series of the Mathematics
Department} "V. Volterra" of the University of Ancona, Via delle Brecce Bianche, Ancona, ITALY (1996).

\bibitem{PPCM02} V. Pravda, A. Pravdov\'a, A. Coley, and R. Milson,
All spacetimes with vanishing curvature invariants, {\it
Classical Quantum Gravity} {\bf 19} (2002), 6213--6236.



\bibitem{S60}I. M. Singer, Infinitesimally homogeneous spaces,
{\it Commun. Pure Appl. Math.} {\bf 13} (1960), 685--697.

\bibitem{W46} H. Weyl, {\bf The Classical Groups} Princeton University Press, 
Princeton, 1946.
\end{thebibliography}
\end{document}